\documentclass[a4paper]{article}

\usepackage[utf8]{inputenc}
\usepackage[T1]{fontenc}
\usepackage{amsmath, amssymb, amsthm}
\usepackage{mathtools}
\usepackage{mathdots}
\usepackage{paralist}
\usepackage[margin=.9in]{geometry}
\usepackage[toc,page]{appendix}
\usepackage{algorithm,algorithmic}

\usepackage{graphicx}
\usepackage[shortlabels]{enumitem}
\usepackage{xcolor}
\usepackage{hyperref}
\usepackage[capitalise,nameinlink]{cleveref}
\crefname{lemma}{Lemma}{Lemma}
\crefname{figure}{Figure}{Figure}
\usepackage{subcaption}

\colorlet{red}{red!80!black}
\colorlet{blue}{blue!70!black}
\colorlet{green}{green!50!black}

\newtheorem{theorem}{Theorem}[section]
\newtheorem{lemma}[theorem]{Lemma}
\newtheorem{remark}[theorem]{Remark}
\newtheorem{example}[theorem]{Example}

\newcommand{\C}{\mathbb{C}}
\newcommand{\cH}{\mathcal{H}}
\newcommand{\cJ}{\mathcal{J}}
\newcommand{\cP}{\mathcal{P}}
\newcommand{\R}{\mathbb{R}}
\newcommand{\cS}{\mathcal{S}}
\DeclarePairedDelimiter{\norm}{\lVert}{\rVert}
\DeclarePairedDelimiter{\abs}{\lvert}{\rvert}
\DeclareMathOperator{\diag}{diag}
\DeclareMathOperator{\imag}{Im}
\DeclareMathOperator{\real}{Re}
\newcommand{\wt}[1]{\widetilde{#1}}

\newcommand{\bp}{BP}
\newcommand{\bpp}{BP$^+$}
\newcommand{\bpm}{BP$^-$}

\author{Chen Greif\thanks{Department of Computer Science, The University of British Columbia, Vancouver, Canada
  (\texttt{greif@cs.ubc.ca}).}, 
  J\"org Liesen\thanks{Institute of Mathematics, Technische Universit\"at Berlin, Berlin, Germany (\texttt{liesen@math.tu-berlin.de}, \texttt{ramme@math.tu-berlin.de}).}, 
  and Justus Ramme\footnotemark[2]}
\title{On Bramble--Pasciak-type preconditioners and their properties} 
\date{\today}

\begin{document}
\maketitle

\begin{abstract}
Bramble and Pasciak (Math.~Comp., 1988) 
discovered a preconditioning method for linear algebraic systems in saddle point form, which under certain conditions transforms a symmetric indefinite saddle point matrix into a nonsymmetric matrix that is positive definite in some non-standard inner product. Thus, the method ``removes'' the indefiniteness of the problem, while retaining the property that the preconditioned system can be solved using an optimal Krylov subspace method with three-term recurrences. Since this well-known method has some limitations, several variants with potential improvements have been proposed in the literature. In this paper, we attempt a systematic assessment of the mathematical properties and practical efficiency of Bramble--Pasciak-type preconditioners. Our goal is to close some gaps in the understanding of these techniques. We identify and analyze the three main methods in this context, highlighting a variant that offers additional flexibility but has received little attention to date. We perform detailed spectral analyzes of the preconditioned matrices and illustrate our theoretical results on several numerical examples.
\end{abstract}



\section{Introduction} 

We consider linear algebraic systems in saddle point form, 
\begin{equation}\label{eq:saddle-point-sys}
    \cS\begin{bmatrix}x\\y\end{bmatrix}
    =\begin{bmatrix}f\\ g\end{bmatrix},
    \;\;\text{where}\;\;
    \cS=\begin{bmatrix} A& B^T \\B & -C\end{bmatrix},
\end{equation} 
with $A=A^T\in\R^{n\times n}$, $C=C^T\in\R^{m\times m}$, $n> m$, and $B\in\R^{m\times n}$. We assume throughout this paper that $A$ is positive definite ($A\succ0$), $C$ is positive semi-definite ($C\succeq 0$), and $B$ has full (row) rank $m$. Under these assumptions $\cS$ is nonsingular and (highly) indefinite.

Systems of the form \eqref{eq:saddle-point-sys} arise in scientific and engineering applications \cite{BenGolLie2005}. In practical applications $\cS$ is often large and sparse, and the linear algebraic system \eqref{eq:saddle-point-sys} is solved iteratively with Krylov subspace methods. Since $\cS$ is symmetric, there exist optimal (norm-minimizing) methods based on three-term recurrences, 
e.g., MINRES~\cite{PaiSau1975}. 
The indefiniteness of $\cS$, however, often slows down the convergence of such methods; see, e.g.,~\cite[Section~5.7.2]{LieStr2013} for details. 

Numerous different preconditioning methods have been developed to overcome slow convergence of methods such as MINRES when applied to saddle point systems;
see, e.g., \cite[Section~10]{BenGolLie2005} and the references therein for a general classification of techniques and many specific examples. Popular methods include 
indefinite and possibly nonsymmetric preconditioners $\cP$, for which the indefiniteness is ``removed'' in the sense that the eigenvalues of the preconditioned matrix $\cP^{-1}\cS$ are entirely in the right (or left) half plane; see, e.g., \cite{BenSim2006,BraPas1988,Krz2011,LiePar2008,LieRam2024,Pes2011,SchZul2007}.  
Usually, the preconditioned matrix $\cP^{-1}\cS$ is then nonsymmetric, which means that there no longer exist optimal three-term recurrence Krylov subspace methods.
However, in some cases, it is possible to construct an inner product in which the matrix $\cP^{-1}\cS$ is selfadjoint and positive definite. Such an inner product yields an optimal three-term recurrence Krylov subspace method for solving the preconditioned (nonsymmetric) positive definite system. 

The first such method was derived by Bramble and Pasciak~\cite{BraPas1988}.  In order to describe their idea algebraically, we note that
the saddle point matrix $\cS$ can be block-LU factorized  
as
\begin{equation}\label{eqn:block-lu}
\cS=\begin{bmatrix} A& 0 \\B & -S\end{bmatrix}
\begin{bmatrix} I_n & A^{-1}B^T \\0 & I_m\end{bmatrix},
\end{equation}
where the (negative) Schur complement satisfies 
\[
S=C+BA^{-1}B^T\succ0.
\]
Let $A_0=A_0^T\succ0$ and $S_0=S_0^T\succ0$ be approximations of $A$ and $S$ in some sense. Then the Bramble-Pasciak (BP) preconditioner is defined by 
\begin{equation}\label{eq:BP-P}
\cP_{BP}=\begin{bmatrix} A_0& 0 \\B & -S_0\end{bmatrix},
\end{equation}
which may be interpreted as an approximation of the left block lower-triangular matrix in the factorization \eqref{eqn:block-lu}. If $A_0=A$ and $S_0=S$, then $\cP_{BP}^{-1}\cS$ is exactly the block upper-triangular matrix on the right-hand side of \eqref{eqn:block-lu}, and the Bramble-Pasciak method amounts to solving \eqref{eq:saddle-point-sys} by block Gaussian elimination. 

We have 
\[
\cP_{BP}^{-1}=\begin{bmatrix} A_0^{-1}& 0 \\S_0^{-1}BA_0^{-1} & -S_0^{-1}\end{bmatrix},
\]
and the preconditioned matrix is given by
\begin{align}\label{eqn:Sbp}
\cS_{BP}=\cP_{BP}^{-1}\cS&=
\begin{bmatrix}A_0^{-1}A & A_0^{-1}B^T\\ S_0^{-1}BA_0^{-1}(A-A_0) & S_0^{-1}(C+BA_0^{-1}B^T)
\end{bmatrix},
\end{align}
which is in general nonsymmetric. Now assume additionally that $A-A_0\succ0$, i.e., that $A_0$ approximates $A$ ``from below'', then
\begin{equation}\label{eq:BP-H}
\cH_{BP}=\begin{bmatrix} A-A_0& 0 \\0 & S_0\end{bmatrix}
\end{equation}
is symmetric positive definite and hence defines an inner product. The surprising algebraic fact about the Bramble-Pasciak method is that then the matrix
\begin{align*}
\cH_{BP}\cS_{BP}&=
\begin{bmatrix} AA_0^{-1}A-A& (A-A_0)A_0^{-1}B^T \\BA_0^{-1}(A-A_0) & C+BA_0^{-1}B^T\end{bmatrix}\\
&= \begin{bmatrix} I_n& 0 \\BA^{-1} & I_m\end{bmatrix} \begin{bmatrix} AA_0^{-1}A-A&0  \\ 0& C+BA^{-1}B^T \end{bmatrix}\begin{bmatrix} I_n & A^{-1}B^T \\& I_m\end{bmatrix}
\end{align*}
is also symmetric positive definite. Consequently, 
there exists an optimal three-term recurrence (CG-like) algorithm in the $\cH_{BP}$-inner product for solving the preconditioned linear algebraic system. This is the original method of Bramble and Pasciak. We point out that their original paper considers $S_0=I$ and does not present the method in terms of block preconditioning. The latter was first done, to the best of our knowledge, by Klawonn in \cite{Kla1998}, who considered a general $S_0=S_0^T\succ0$ and also proved spectral bounds for the preconditioned matrices. 

A key challenge with the Bramble-Pasciak method is to find a matrix $A_0$ such that $A-A_0\succ0$. 
Mathematically, as long as $A\succ0$, 
this can be easily achieved by trivial choices such as $A_0=\frac12 A$. 
But if $A$ is ill conditioned, then so is $\cH_{BP}$.  Even the most sophisticated choice of $A_0$ can lead to an ill-conditioned inner product matrix $\cH_{BP}$. 
This situation has motivated several modifications and extensions of the original method of Bramble and Pasciak. For example, Stoll and Wathen \cite{StoWat2008} proposed 
the {\em Bramble-Pasciak-plus} (\bpp{}) preconditioner, which uses $A+A_0$ instead of $A-A_0$ in the inner product matrix, as well as a corresponding modification in the preconditioner itself. One advantage of this technique is, of course, that $A+A_0\succ0$ whenever $A_0\succ0$. Other modifications of the Bramble-Pasciak method were given, e.g., in \cite{Pes2011,SchZul2007}.

The goal of this paper is to close some gaps in the understanding of the mathematical properties and practical efficiency of Bramble--Pasciak-type preconditioners. In a nutshell, and disregarding trivial modifications through additional parameter choices, the three main methods in this context can be characterized by choosing the (1,1) block of the inner product matrix as $A-A_0$ (original BP), $A+A_0$ (\bpp{}), or $A_0-A$ (\bpm{}), with the corresponding suitable modifications of the preconditioner. We analyze these three variants, and this allows us to work out their main features and differences, which we then illustrate in numerical experiments. We particularly highlight the \bpm{} method, which gives extra flexibility in situations where $A$ is ill-conditioned (or even singular) and which has not been studied in the literature so far.

The three methods we study in this paper are in principle contained in the general parameterized approaches in~\cite{DolGouStoWat2010,Krz2011}. Moreover, variants of some of the eigenvalue bounds we present for the original BP and the \bpp{} methods can be found in the literature (references are given in the appropriate places). However, the paper contains several new results, in particular for the \bpm{} method, and the first systematic analysis and comparison of the main characteristics of Bramble-Pasciak-like methods we are aware of.

The paper is organized as follows. In Section \ref{sec:BP-variants} we provide the necessary background from the theory of Krylov subspace methods on the existence of optimal three-term recurrence methods, and we derive the three main Bramble-Pasciak-like methods from a general framework. We then show the rather surprising result that the preconditioned matrix in the \bpm{} method is diagonalizable with real eigenvalues under very mild conditions. Section \ref{sec:BP-analysis} contains a detailed spectral analysis of the preconditioned matrices for the three methods. The essential differences are first worked out on a simple one-parameter case. This is followed by a general spectral analysis that includes studying the tightness of the eigenvalue bounds and the convergence behavior of the optimal three-term recurrence methods defined by the different preconditioning techniques. 
In Section \ref{sec:examples} we illustrate our theoretical results on several numerical examples. Finally, in Section \ref{sec:conclusion}, we summarize our findings and draw some conclusions.

\section{Variants of the Bramble-Pasciak preconditioner} \label{sec:BP-variants}

We start with the following well known theorem, which forms the theoretical foundation of the existence of optimal (in the sense of norm-minimizing) three-term recurrence Krylov subspace methods for a given matrix in a prescribed inner product. The original result goes back to Faber and Manteuffel \cite{FabMan1984}. For the formulation below and, in particular, a precise definition of the Krylov subspace methods considered, we refer to \cite{LieStr2008} or \cite[Chapter 4]{LieStr2013}. 

\begin{theorem}\label{thm:optimal}
    Let $M\in\C^{n\times n}$ and let $H\in\C^{n\times n}$ be Hermitian positive definite. Then the following are equivalent:
    \begin{enumerate}[(i)]
    \item An optimal three-term recurrence Krylov subspace method exists for $M$ in the $H$-inner product.
	\item $M$ is $H$-normal(1), i.e., the $H$-adjoint of $M$ satisfies $M^+=H^{-1}M^H H=p(M)$ for a polynomial $p$ of degree one.
	\item $M$ is diagonalizable, $M=W\Lambda W^{-1}$ (where without loss of generality we assume that the eigenvalues and eigenvectors are ordered so that equal eigenvalues form a single diagonal block in $\Lambda$) with collinear eigenvalues, and $H^{-1}=WDW^H$, where $D$ is Hermitian positive definite and block-diagonal with block sizes as those in $\Lambda$.
    \end{enumerate}
\end{theorem}

In the terminology of Theorem~\ref{thm:optimal}, if the matrix $\cH_{BP}$ in \eqref{eq:BP-H} is positive definite, then the Bramble-Pasciak preconditioned matrix $\cS_{BP}$ is $\cH_{BP}$-normal(1). This fact was already explicitly pointed out by Klawonn in \cite[p. 177]{Kla1998}. \Cref{sec:methods} contains two important examples of optimal three-term recurrence methods for general $H$-selfadjoint matrices (i.e., when $M^+=M$), namely the $H$-PCG and the $H$-PMINRES method. The former is applicable when $M$ is positive definite in the $H$-inner product, while the latter only requires that $M$ be $H$-selfadjoint.

Theorem~\ref{thm:optimal} motivates the question whether for a given (real) saddle point matrix $\cS$ there exists a preconditioner $\cP$ and a symmetric positive definite (i.e., inner product) matrix $\cH$ so that $\cP^{-1}\cS$ is $\cH$-normal(1). Then the preconditioned matrix is diagonalizable with collinear eigenvalues. Ideally, these eigenvalues are all real and positive (or all real and negative), so that the indefiniteness of the problem is ``removed.'' 

The following result, which is nothing more than a parameterized form of the original Bramble-Pasciak method, is our starting point for studying this question. 
\begin{lemma}\label[lemma]{lem:basic}
    Let $\cS$ be as in \eqref{eq:saddle-point-sys} with the corresponding assumptions.
    Let $A_0=A_0^T\in\R^{n\times n}$ and $S_0=S_0^T\in\R^{m\times m}$ be nonsingular, $\alpha,\beta\in\R$, $\beta\neq 0$, and define
    \begin{equation}\label{eqn:PandH}
    \cP_{\alpha,\beta}=\begin{bmatrix} A_0& 0 \\\alpha B & \beta S_0\end{bmatrix} \quad\mbox{and}\quad
    \cH_{\alpha,\beta}=\begin{bmatrix} \frac{1}{\beta} (A_0-\alpha A)& 0 \\0 & S_0\end{bmatrix}.
    \end{equation}
    Then the following holds:
    \begin{enumerate}[(i)]
    \item The preconditioner $\cP_{\alpha,\beta}$ is nonsingular, and the preconditioned matrix  
    \begin{equation}
    \begin{split}
    \cS_{\alpha,\beta} &=
    \cP_{\alpha,\beta}^{-1}\cS=
    \begin{bmatrix}
        A_0^{-1} & 0\\ -\frac{\alpha}{\beta}S_0^{-1}BA_0^{-1} & \frac{1}{\beta}S_0^{-1}
    \end{bmatrix}
    \begin{bmatrix} A & B^T\\B & -C \end{bmatrix} \\
    &=
    \begin{bmatrix} A_0^{-1}A & A_0^{-1}B^T \\
    \frac{1}{\beta}S_0^{-1}B(I_n-\alpha A_0^{-1}A) & -\frac{1}{\beta}S_0^{-1}(C+\alpha BA_0^{-1}B^T)\end{bmatrix} \label{eqn:PinvS}
    \end{split}
    \end{equation}
    satisfies
    \begin{align*}
    \cH_{\alpha,\beta}\cS_{\alpha,\beta}&=
    \frac{1}{\beta} \begin{bmatrix} A-\alpha AA_0^{-1}A & (I_n-\alpha AA_0^{-1})B^T \\B(I_n-\alpha A_0^{-1}A) & -(C + \alpha BA_0^{-1}B^T)\end{bmatrix}. 
    \end{align*}
    In particular, $\cH_{\alpha,\beta}\cS_{\alpha,\beta}$ is symmetric, i.e., the preconditioned matrix $\cS_{\alpha,\beta}$ is $\cH_{\alpha,\beta}$-selfadjoint.  
    \item
    If $\cH_{\alpha,\beta}\succ0$, i.e., if
    $\frac{1}{\beta} (A_0-\alpha A)\succ0$ and $S_0\succ0$,
    then the preconditioned matrix $\cS_{\alpha,\beta}$ is diagonalizable with only real eigenvalues (though not necessarily all positive or all negative).
    \end{enumerate}
\end{lemma}

\begin{proof}~
\begin{enumerate}[(i)]
    \item 
    The calculation showing that $\cH_{\alpha,\beta}\cS_{\alpha,\beta}$ has the required form is straightforward.

    \item 
     The symmetry of $\cH_{\alpha,\beta}\cS_{\alpha,\beta}$ is equivalent to $\cH_{\alpha,\beta} \cS_{\alpha,\beta} =\cS_{\alpha,\beta}^T\cH_{\alpha,\beta}$. Thus, if $\cH_{\alpha,\beta}\succ0$, then 
     \[\cS_{\alpha,\beta}^+=\cH_{\alpha,\beta}^{-1}\cS_{\alpha,\beta}^H\cH_{\alpha,\beta}=\cH_{\alpha,\beta}^{-1}\cS_{\alpha,\beta}^T\cH_{\alpha,\beta}=\cS_{\alpha,\beta},
     \]
    i.e., $\cS_{\alpha,\beta}$ is $\cH_{\alpha,\beta}$-normal(1), and hence diagonalizable with collinear eigenvalues by Theorem~\ref{thm:optimal}. 
    Since $\cS_{\alpha,\beta}$ is $\cH_{\alpha,\beta}$-selfadjoint,
    the eigenvalues of $\cS_{\alpha,\beta}$ are indeed all real. 
\end{enumerate}
\end{proof}
We note that statement (i) in \cref{lem:basic} holds for general $A=A^T\in\R^{n\times n}$ and $C=C^T\in\R^{m\times m}$.
\begin{remark}
The choice $\alpha=0$ and $\beta=1$ in Lemma~\ref{lem:basic} yields the popular block-diagonal preconditioner
\begin{equation*}
    \cP_{0,1} = 
    \begin{bmatrix}
        A_0 & 0\\ 0 & S_0
    \end{bmatrix};
\end{equation*}
 see~\cite[Section~10.1]{BenGolLie2005} for general theory. If $A_0,S_0\succ0$, then the corresponding preconditioned matrix $\cS_{0,1}$ is diagonalizable with real eigenvalues but still indefinite. If $C=0$, the choice $A_0=A$ and $S_0=S$ yields the well known ``ideal'' block-diagonal or Murphy-Golub-Wathen~\cite{MurGolWat2000} preconditioner with $\sigma(\cS_{0,1})=\{(1-\sqrt{5})/2,1,(1+\sqrt{5})/2\}$. We will use block-diagonal preconditioners as a base case in the numerical experiments.
\end{remark}
Suppose that $A_0\succ0$ and $\cH_{\alpha,\beta}\succ0$, i.e., $\frac{1}{\beta} (A_0-\alpha A)\succ0$ and $S_0\succ0$, and factorize 
\begin{multline}\label{eqn:fact1}
    \cH_{\alpha,\beta} \cS_{\alpha,\beta}
    \\ =
    \begin{bmatrix}
        I_n & 0 \\ BA^{-1} & I_m 
    \end{bmatrix}
    \begin{bmatrix}
        \frac{1}{\beta}A(A^{-1}-\alpha A_0^{-1})A & 0 \\ 0 & -\frac{1}{\beta}(C+BA^{-1}B^T)
    \end{bmatrix}
    \begin{bmatrix}
        I_n & A^{-1}B^T \\ 0 & I_m 
    \end{bmatrix}.
\end{multline}
Since $A,A_0\succ0$, we have that  $\frac{1}{\beta} A(A^{-1} - \alpha A_0^{-1})A\succ0$ is equivalent to $\frac{1}{\beta} (A_0-\alpha A)\succ0$, and thus the (1,1)-block of the block-diagonal matrix is always positive definite. Thus, the question whether $\cS_{\alpha,\beta}$ is positive definite or indefinite in the $\cH_{\alpha,\beta}$ inner product is completely determined by the (2,2)-block, and hence by the sign of $\beta$. Choosing $\beta<0$ gives $\cH_{\alpha,\beta}\cS_{\alpha,\beta}\succ0$, so $\cS_{\alpha,\beta}$ is diagonalizable with positive eigenvalues.

Three choices of $(\alpha,\beta)$ are of particular interest; we list them in Table~\ref{tab:precs}. The combination $\alpha=\beta=-1$ is not suitable, because the resulting (1,1)-block of the matrix $\cH_{-1,-1}$ is given by $-(A_0+A)$, which is negative definite. Hence, in this case $\cH_{-1,-1}$ is indefinite and does not define an inner product.

\begin{table}[h]
\begin{center}
\begin{tabular}{c|c|l}
        $(\alpha,\beta)$ & (1,1)-block of $\cH_{\alpha,\beta}$ &  Method \\ \hline
        $(1,-1)$ & $A-A_0$ & Bramble-Pasciak (\bp{}) \cite{BraPas1988}  \\
        $(-1,1)$ & $A+A_0$ & Bramble-Pasciak-plus (\bpp{}) \cite{StoWat2008} \\
        $(1,1)$ & $A_0-A$ & Bramble-Pasciak-minus (\bpm{})  
\end{tabular}  
\end{center}
\caption{Three BP-type preconditioners of interest.}\label{tab:precs}
\end{table}

Note that the preconditioners $\cP_{\alpha,\beta}$ resulting from the choices in Table~\ref{tab:precs} are related by
\begin{equation*}
\cP_{1,-1}=\cJ \cP_{-1,1}=\cP_{1,1}\cJ,\quad \mbox{where}\quad
\cJ=\begin{bmatrix} I_n & 0\\0 &-I_m\end{bmatrix}.
\end{equation*}
As mentioned in the Introduction, all 
three methods, as well as several other block preconditioners and the corresponding inner product matrices (possibly indefinite), are contained in the parameterization given in~\cite{Krz2011}. In particular, the method we call here ${\rm BP}^{-}$ is obtained by setting $\delta=c=1$ and $d=0$ in that parameterization. However, the paper~\cite{Krz2011} does not contain detailed spectral analyses of the corresponding preconditioned matrices.

As shown in \cref{lem:basic}, when the matrix $\cH_{\alpha,\beta}$ is positive definite  the eigenvalues of the preconditioned matrix $\cS_{\alpha,\beta}$ are all real.
As discussed in the Introduction, in the case of the original BP method the eigenvalues are 
also all positive, making the matrix $\cS_{1,-1}$ positive definite in the $\cH_{1,-1}$-inner product.  
In contrast, the matrix $\cS_{-1,1}$ in the \bpp{} method is indefinite in the $\cH_{-1,1}$-inner product, which was already noted in \cite{StoWat2008}.

Interestingly, the matrix $\cS_{1,1}$ in the \bpm{} method has a real spectrum and is diagonalizable even if $\cH_{1,1}$ is not positive definite (i.e., when the condition for it being an inner product  is violated). The realness of the spectrum was first observed in \cite[Theorem~2.1]{Cao2007}. In the following, we also prove diagonalizability, and note that our proof does not use the fact that $B$ is of full rank.

\begin{theorem}\label{thm:BPm_real_eig}
Let $\cS$ be as in \eqref{eq:saddle-point-sys} with the corresponding assumptions. In the notation of Lemma~\ref{lem:basic}, assume that $A_0,S_0\succ0$.     Then the \bpm{} preconditioned matrix $\cS_{1,1}$ is diagonalizable with real eigenvalues (even when $\cH_{1,1}$ is indefinite, i.e., $A_0-A$
is not positive definite).
\end{theorem}

\begin{proof}
    We start by showing that the eigenvalues of $\cS_{1,1}$ are real.
    In order to do so, note that the eigenvalues of the matrix $\cS_{1,1}$ are the same as the ones of the pencil 
    \begin{equation*}
        \cS 
        \begin{bmatrix}
            x\\y
        \end{bmatrix}
        =
        \lambda
        \cP_{1,1}
        \begin{bmatrix}
            x\\y
        \end{bmatrix}.
    \end{equation*}
    We thus obtain the two equations 
    \begin{equation}\label{eq:BPm-eig-eqs}
    \begin{aligned}
        B^Ty &= (\lambda A_0-A)x \\  
        (1-\lambda)Bx &= (\lambda S_0+C)y.
    \end{aligned}
    \end{equation}
    Multiplying the first equation on the left by $x^H$, the second by $y^H$ and substituting the Hermitian transpose of the first equation into the second  yields 
    \begin{equation}\label{eq:BPm-combined-eig-eq}
    \begin{aligned}
        0 &= (1-\lambda)x^H(\overline{\lambda} A_0-A)x - y^H(\lambda S_0 +C)y
        \\
        &= - \abs{\lambda}^2 x^HA_0x + \lambda (x^HAx - y^HS_0y ) + \overline{\lambda}x^HA_0x -x^HAx - y^HCy. 
    \end{aligned}
    \end{equation}
    Taking the imaginary part yields
    \begin{equation}\label{eq:BPm-imag-eig}
        0 = \imag(\lambda)(x^HAx - y^HS_0y - x^HA_0x). 
    \end{equation}
    Thus, if we assume that $\imag(\lambda)\neq 0$, we must have $0 = x^H(A_0-A)x + y^HS_0y$. Taking the real part of \eqref{eq:BPm-combined-eig-eq} and using \eqref{eq:BPm-imag-eig}, we get 
    \begin{equation}\label{eq:BPm-H-neutral-eig-vec}
    \begin{aligned}
        0 &= -\abs{\lambda}^2 x^HA_0x + \real(\lambda)(x^HAx-y^HS_0y + x^HA_0x) -x^HAx -y^HCy 
        \\
        &= -((\real(\lambda)-1)^2 +\imag(\lambda)^2 )x^HA_0x - y^HS_0y -y^HCy.
    \end{aligned}
    \end{equation}
    Since $A_0\succ0$ and $S_0\succ0$ and we assumed that $\imag(\lambda)\neq0$, it follows that $x=0$ and $y=0$, which is a contradiction to $[x^H\;y^H]^H$ being an eigenvector. Thus $\imag(\lambda ) = 0$. 
    
    To prove the diagonalizability, we proceed as follows. Let $\lambda$ be an eigenvalue of $\cS_{1,1}$ and  $w = [u^T\;v^T]^T\in \ker((\lambda I_{n+m}-\cS_{1,1})^2)$. Then $z= [x^T\;y^T]^T =(\lambda I_{n+m}-\cS_{1,1})w\in \ker (\lambda I_{n+m}-\cS_{1,1})$ and we have  
    \begin{equation*}
        x^T(A_0-A)x + y^TS_0y
        = z^T\cH_{1,1}z 
        = w^T\cH_{1,1}(\lambda I_{n+m}-\cS_{1,1})^2w = 0.
    \end{equation*}
    Note that this is exactly the condition we derived earlier from \eqref{eq:BPm-imag-eig} for $\imag(\lambda)\neq0$ and hence, we similarly obtain \eqref{eq:BPm-H-neutral-eig-vec}, however, this time for real $\lambda$, i.e.,
    \begin{equation*}
        0=-(\lambda-1)^2x^TA_0x -y^TS_0y -y^TCy.
    \end{equation*}
    If $\lambda\neq1$, then $z=[x^T\;y^T]^T=0$ and hence $w$ is already an eigenvector of $\cS_{1,1}$. If, on the other hand, $\lambda=1$, then we only obtain $y=0$. We can see from \eqref{eq:BPm-eig-eqs} that in this case $A_0x=Ax$. In order to further show that $x=0$, we use  $z=( I_{n+m}-\cS_{1,1})w$, which yields the two equations 
    \begin{subequations}
    \label{eq:S11-eig-eqs}
    \begin{align}
        x &= (I_n - A_0^{-1}A)u -A_0^{-1}B^Tv , \label{eq:S11first}\\ 
        0 &= -S_0^{-1}B(I_n-A_0^{-1}A)u + v +S_0^{-1}(BA_0^{-1}B^T+C)v.
        \label{eq:S11second}
        \end{align}
    \end{subequations}
    Rearranging and substituting \eqref{eq:S11first} into \eqref{eq:S11second} for $(I_n-A_0^{-1}A)u$ gives 
\begin{equation}\label{eq:S11-combined-eig-eq}
        0 = -S_0^{-1}B(x+A_0^{-1}B^Tv) + v +S_0^{-1}(BA_0^{-1}B^T+C)v 
        = -S_0^{-1}Bx+  v +S_0^{-1}Cv. 
    \end{equation}
    Multiplying \eqref{eq:S11first} from the left by $x^TA_0$ and \eqref{eq:S11-combined-eig-eq}  by $v^TS_0$ yields 
    \begin{equation*}
    \begin{aligned}
        x^TA_0x &= -x^TB^Tv,
        \\
        0 &= -v^TBx + v^TS_0v + v^TCv,
    \end{aligned}
    \end{equation*}
    which gives us
    \begin{equation*}
        0 = x^TA_0x + v^TS_0v + v^TCv.
    \end{equation*}
    Thus, we again obtain $z=0$ and  $w$ is an eigenvector of $\cS_{1,1}$. Therefore, $\cS_{1,1}$ is diagonalizable as required.
\end{proof}
\begin{remark} \cref{thm:BPm_real_eig} 
shows that the \bpm{} preconditioned matrix $\cS_{1,1}$ is diagonalizable with real eigenvalues even under weak conditions, which in theory allows an optimal three-term recurrence method by \cref{thm:optimal}. 
However, there does not seem to be an easy way to construct the corresponding inner product, so this method can hardly be realized in practice.  
\end{remark}

\section{Spectral analysis of the preconditioned matrices} \label{sec:BP-analysis}

This section contains a spectral analysis of the preconditioned matrices corresponding to 
the three different BP-type methods given in Table \ref{tab:precs}. We first present a simple analysis where $C=0$ and the three methods are controlled by a single parameter through $A_0=\gamma A$ with $\gamma>0$. This reveals some of the main effects and sets the stage for a more general analysis that follows in Section~\ref{sec:BP-general-analysis} and which contains fairly general approximations $A_0\approx A$ and $S_0\approx S$. The analysis is based, to a large extent, on eigenvalue differences of $A$ and $A_0$ as well as $S$ and $S_0$. In practice, small differences are guaranteed by \emph{spectrally equivalent approximations}. Recall that two symmetric matrices $M_1,M_2\in\R^{n\times n}$ are spectrally equivalent when there exist constants $\delta,\Delta>0$ such that 
\begin{equation*}
    \delta \cdot u^TM_2u \le u^TM_1u \le \Delta \cdot u^TM_2u,
\end{equation*}
for all  $ u\in\R^{n}$.

\subsection{One-parameter analysis of a special case}\label{sec:special-case}

If $A_0=A$ and $S_0=S$, then the original BP method as well as the \bpm{} method become block Gaussian elimination, since in this case 
\begin{equation*}
    \cS_{1,-1} = 
    \begin{bmatrix}
        I_n & A^{-1}B^T \\ 0 & I_m
    \end{bmatrix}
\quad \text{and} \quad
\cS_{1,1} = 
    \begin{bmatrix}
        I_n & A^{-1}B^T \\ 0 & -I_m
    \end{bmatrix}.
\end{equation*}
In this sense we may consider these choices as ``ideal preconditioners.'' Note, however, that for $A_0=A$ the matrices $\cH_{1,-1}$ (${\rm BP}$) and $\cH_{1,1}$ (${\rm BP}^{-}$) both have a zero (1,1)-block. For a usable technique that potentially leads to an optimal three-term recurrence method we will therefore consider a scaled variant, $A_0=\gamma A$, for some parameter $\gamma>0$ .

\begin{lemma}\label[lemma]{lem:ideal_prec}
    Let $\cS$ be as in \eqref{eq:saddle-point-sys} with the corresponding assumptions, and additionally let $C=0$. In the notation of Lemma~\ref{lem:basic}, assume that $A_0=\gamma A\succ0$ and $S_0\succ0$. Let the preconditioned matrix $\cS_{\alpha,\beta}$ be as in \eqref{eqn:PinvS}. In addition, let $s_1,\dots,s_m$ be the $m$ (real and positive) eigenvalues of $S_0^{-1}BA^{-1}B^T$.
    Then 
    \begin{equation}\label{eqn:nu-pm}
    \begin{aligned}
    \sigma(\cS_{\alpha,\beta}) =& 
    \left\{\frac{1}{\gamma}\right\} \cup
    \left\{
    \frac{\nu_{\pm}^{(j)}}{\gamma} 
    \bigg\vert\; j = 1,\dots,m
    \right\}
    \quad \text{with}
    \\
    \nu_{\pm}^{(j)} =& 
    \frac{\beta-\alpha s_j}{2\beta} \pm \sqrt{\frac{(\beta-\alpha s_j)^2}{4\beta^2}+\frac{\gamma s_j}{\beta}}.
    \end{aligned} 
    \end{equation}
    Furthermore, $\cS_{\alpha,\beta}$ is diagonalizable when $\gamma \neq - \frac{(\beta-\alpha s_j)^2}{4\beta s_j}$ for all $j=1,\dots,m$. 
\end{lemma}
\begin{proof}
We immediately see that 
\begin{equation*}
\gamma\cS_{\alpha,\beta} = 
\begin{bmatrix}
I_n &  A^{-1}B^T \\ \frac{\gamma-\alpha}{\beta}S_0^{-1}B & -\frac{\alpha}{\beta} S_0^{-1}BA^{-1}B^T
\end{bmatrix}.
\end{equation*}
Let $V\in \R^{n\times n-m}$ be an $A$-orthonormal basis of $\ker(B)$. Then the matrix 
\begin{equation*}
    X_1 = 
    \begin{bmatrix}
        V & A^{-1}B^T & 0 \\ 0 & 0 & I_m
    \end{bmatrix}
\end{equation*}
is nonsingular, since 
\begin{equation}
    X_1^T \diag(A,I_m) X_1 = \diag (I_{n-m}, BA^{-1}B^T, I_m),
\end{equation}
hence $X_1$ is nonsingular. A calculation reveals that 
\begin{equation*}
    \gamma X_1^{-1}S_{\alpha,\beta}X_1 = 
    \begin{bmatrix}
        I_{n-m} & 0 & 0 \\ 0 & I_m & I_m \\ 0 & \frac{\gamma-\alpha}{\beta}S_0^{-1}BA^{-1}B^T & -\frac{\alpha}{\beta} S_0^{-1}BA^{-1}B^T
    \end{bmatrix}.
\end{equation*}
Let now $U\in\R^{m\times m}$ be an eigenvector basis of $S_0^{-1}BA^{-1}B^T$, such that 
\[
U^{-1}S_0^{-1}BA^{-1}B^TU = D_S = \diag(s_1,\dots,s_m).
\] 
Setting $X_2=\diag(I_{n-m}, U,U)$, we have
\begin{equation*}
    \gamma (X_1X_2)^{-1}\cS_{\alpha,\beta}X_1X_2 = 
    \begin{bmatrix}
        I_n & 0 & 0 \\ 0 & I_m & I_m \\ 0 & \frac{\gamma-\alpha}{\beta}D_S & -\frac{\alpha}{\beta} D_S
    \end{bmatrix}.
\end{equation*}
Thus, $\cS_{\alpha,\beta}$ has an eigenvalue $\frac{1}{\gamma}$ with multiplicity $n-m$ and the remaining eigenvalues are given by the eigenvalues of the matrices 
\begin{equation*}
\wt{\cS}_j=
\frac{1}{\gamma}
    \begin{bmatrix}
        1 & 1 \\  \frac{\gamma-\alpha}{\beta}s_j & -\frac{\alpha}{\beta} s_j
    \end{bmatrix},
\end{equation*}
which are exactly of the form  \eqref{eqn:nu-pm}.
Furthermore, $\nu_+^{(j)}\neq \nu_-^{(j)}$ if $\gamma \neq - \frac{(\beta-\alpha s_j)^2}{4\beta s_j}$ and, hence, $\wt{\cS}_j$ is diagonalizable in this case. If, on the other hand, $\gamma = - \frac{(\beta-\alpha s_j)^2}{4\beta s_j}$ then $\wt{\cS}_j$ has a duplicate eigenvalue, but is not a multiple of the identity matrix and therefore cannot be diagonalizable. 
In particular, $\cS_{\alpha,\beta}$ is diagonalizable if and only if $\wt{\cS}_j$ is diagonalizable for all $j=1,\dots,m$.
\end{proof}

Note that for any $\beta>0$ we have
$$-\frac{(\beta-\alpha s_j)^2}{4\beta s_j}\leq 0.$$
Hence, for any $\gamma>0$  
the matrix $\cS_{\alpha,\beta}$ is diagonalizable.
In view of \eqref{eqn:fact1}, the choice $\beta>0$ 
corresponds to the case where the matrix $\cS_{\alpha,\beta}$ is (still) indefinite and when $\cS$ and $\cS_{\alpha,\beta}$ have the same inertia.

Applying Lemma~\ref{lem:ideal_prec} to the three BP-type methods with the exact Schur complement $S_0 = BA^{-1}B^T$ (i.e., $s_j=1$ for $j=1,\dots,m$ in \cref{lem:ideal_prec}) yields
\begin{equation}\label{eq:eigvals}
\begin{aligned}
    {\rm BP}:&&\sigma(\cS_{1,-1}) & = \left\{ \frac{1}{\gamma},\quad 
    \frac{1}{\gamma}(1 + \sqrt{1-\gamma}),\quad 
    \frac{1}{\gamma} (1 - \sqrt{1-\gamma})\right\},
    \\ 
    {\rm BP}^{+}:&& \sigma(\cS_{-1,1}) & = \left\{ \frac{1}{\gamma},\quad 
        \frac{1}{\gamma} (1 + \sqrt{1+\gamma}),\quad 
        \frac{1}{\gamma} (1 - \sqrt{1+\gamma})  
        \right\},
    \\ 
    {\rm BP}^{-}:&& \sigma(\cS_{1,1}) & =\left\{ \frac{1}{\gamma},\quad 
        \frac{1}{\sqrt{\gamma}},\quad 
        - \frac{1}{\sqrt{\gamma}}  
        \right\}.
\end{aligned}
\end{equation}
In the original BP method we require $A-A_0=(1-\gamma)A\succ0$ and hence $\gamma\in (0,1)$. If $\gamma\uparrow 1$, then $\sigma(\cS_{1,-1})$ approaches $\{1\}$. In the \bpp{} method we do not have restrictions on $\gamma>0$, and for $\gamma\approx 1$ we have $\sigma(\cS_{-1,1})$ close to $\{1-\sqrt{2},1, 1+\sqrt{2}\}$, which structurally is close to $\sigma(\cS_{0,1}) = \{(1-\sqrt{5})/2,1, (1+\sqrt{5})/2\}$; cf.~Section~\ref{sec:BP-variants}. Finally, for the \bpm{} method we require $A_0-A=(\gamma-1)A\succ0$ and hence $\gamma>1$. If $\gamma\downarrow 1$, then $\sigma(\cS_{1,1})$ approaches $\{-1,1\}$. In Example~\ref{ex:exact_stokes} we will illustrate and discuss how the behavior of the optimal iterative solvers corresponds to these observations about the eigenvalue distributions of the different BP-type preconditioned matrices.

Note that the values $\nu_{\pm}^{(j)}$ for the \bpm{} method, where $\alpha=\beta=1$, are given by
\begin{align*}
    \nu_+^{(j)} &= 1 + \frac{1+s_j}{2}\left(\sqrt{1+(\gamma-1)\frac{4s_j}{(1+s_j)^2}}-1\right),\\
    \nu_-^{(j)} &= -s_j - \frac{1+s_j}{2}\left(\sqrt{1+(\gamma-1)\frac{4s_j}{(1+s_j)^2}}-1\right). 
\end{align*}
These expressions illustrate the following observation, which we will generalize in Section~\ref{sec:BP-general-analysis} and numerically illustrate in~\cref{ex:exact_stokes}: If $A_0$ is a close approximation of $A$, in this special case if $\gamma\approx 1$, then $\gamma-1\approx 0$ and the term in brackets in $\nu_{+}^{(j)}$ is close to $0$. Thus, as long as the $s_j$ are of moderate size, we will have $\nu_{+}^{(j)}\approx 1$, and hence the positive eigenvalues of the \bpm{} preconditioned matrix will be close to $1$. However, under the same assumptions we expect a wider spread of the negative eigenvalues since $\nu_{-}^{(j)}\approx -s_j$, which is illustrated in Figure~\ref{fig:exact_stokes_eig_BP} (c).

In contrast, for the \bpp{} method, where $\alpha=-1$ and $\beta=1$, we have
\begin{align*}
    \nu_+^{(j)} &= \frac{1 + s_j}{2} \left(1 +\sqrt{1+\gamma \frac{4s_j}{(1+s_j)^2}}\right),\\
    \nu_-^{(j)} &= \frac{1 + s_j}{2} \left(1 -\sqrt{1+\gamma \frac{4s_j}{(1+s_j)^2}}\right).
\end{align*}
Here both quantities are significantly influenced by the sizes of the $s_j$, and we expect a wider spread of the negative and the positive eigenvalues even when $A_0\approx A$ (here $\gamma\approx 1$), which is illustrated in Figure~\ref{fig:exact_stokes_eig_BP}~(b).

\subsection{General spectral analysis}\label{sec:BP-general-analysis}

Eigenvalue bounds for the preconditioned matrix in the original Bramble-Pasciak method (the matrix $\cS_{1,-1}$), in some cases even for an indefinite matrix $A-A_0$, were given in~\cite{Not2014,Sim2004,Zul2002}, with the most general result given by~\cite[Theorem 4.3]{Not2014}. For comparison with our bounds for the preconditioned matrices in the \bpp{} and \bpm{} methods below, we state the following result, which is an immediate consequence of~\cite[Lemma~5]{Sim2004}. 

\begin{theorem}
\label{thm:eig_BP}
Let $\cS$ be as in \eqref{eq:saddle-point-sys} with the corresponding assumptions. In the notation of Lemma~\ref{lem:basic}, assume that $A_0,S_0,\cH_{1,-1}\succ0$.
Then the eigenvalues of the \bp{} preconditioned matrix
$\cS_{1,-1}$ lie in the interval $[a,b]\subseteq\R_+$ with 
    \begin{align*}
        a &= \frac{\min\{1,\lambda_{\min}(S_0^{-1}(C+BA_0^{-1}B^T))\}}{2\lambda_{\max}(A_0^{-1}A)(1+\lambda_{\max}(S_0^{-1}(C+BA_0^{-1}B^T)))} \quad\mbox{and}
        \\ 
        b &= 2\max \{\lambda_{\max}(A_0^{-1}A),\lambda_{\max}(S_0^{-1}(C+BA_0^{-1}B^T))\}.
    \end{align*}
\end{theorem}
\begin{proof}
    Since $\cH_{1,-1}\succ0$ we have $A-A_0\succ0$ and therefore $A_0^{-\frac{1}{2}}AA_0^{-\frac{1}{2}}-I_n\succ0$, where we have used  $A_0\succ0$. Thus, the matrix 
    $X = \diag((A_0^{-\frac{1}{2}}AA_0^{-\frac{1}{2}}-I_n)^{\frac12}A_0^{\frac{1}{2}}, S_0^{\frac12})$ is well defined and yields 
    \begin{equation*}
        X\cS_{1,-1}X^{-1} = 
        \begin{bmatrix}
            A_0^{-\frac{1}{2}}AA_0^{-\frac{1}{2}} &
            (A_0^{-\frac{1}{2}}AA_0^{-\frac{1}{2}}-I_n)^{\frac12}A_0^{-\frac{1}{2}}B^TS_0^{-\frac{1}{2}} \\ 
            S_0^{-\frac{1}{2}}BA_0^{-\frac{1}{2}}(A_0^{-\frac{1}{2}}AA_0^{-\frac{1}{2}}-I_n)^{\frac12} &
            S_0^{-\frac{1}{2}}(C+BA_0^{-1}B^T)S_0^{-\frac{1}{2}}
        \end{bmatrix}.
    \end{equation*}
    The specified eigenvalue bounds can now be read off directly from~\cite[Lemma~5]{Sim2004}.
\end{proof}
As observed in \cite{Sim2004} these bounds indicate that if the matrices $A_0^{-1}A$ and $S_0^{-1}(C+BA_0^{-1}B^T)$ are close to the identity matrix, then the eigenvalues of $\cS_{1,-1}$ will cluster around $1$.

A spectral analysis of the \bpp{} preconditioned matrix (the matrix $\cS_{-1,1}$) was performed in~\cite{StoWat2008}. However, the spectral bounds given there depend on the respective eigenvectors of $\cS_{-1,1}$.
In the following, we will complement the existing analyses with eigenvalue bounds for both $\cS_{-1,1}$ and $\cS_{1,1}$. 
\begin{theorem}\label{thm:eig_bounds}
Let $\cS$ be as in \eqref{eq:saddle-point-sys} with the corresponding assumptions. In the notation of Lemma~\ref{lem:basic}, assume that $A_0,S_0\succ0$. 
Further let 
\begin{align*}
    \lambda^A_{\min} & = \lambda_{\min}(A_0^{-1}A), 
    & 
    \lambda^A_{\max} &= \lambda_{\max}(A_0^{-1}A),
    \\
    \lambda^S_{\min} & = \lambda_{\min}(S_0^{-1}(C+BA_0^{-1}B^T)), 
    & 
    \lambda^S_{\max} &= \lambda_{\max}(S_0^{-1}(C+BA_0^{-1}B^T)),
    \\ 
    \lambda^C_{\max} &= \lambda_{\max}(S_0^{-1}C) 
    & 
    \lambda^T_{\max} &= \lambda_{\max}(A_0^{-1}(A+B^TS_0^{-1}B)),
    \\
    \sigma_{\min} & = \lambda_{\min}(S_0^{-1}(BA_0^{-1}B^T)), 
    & 
    \sigma_{\max} &= \lambda_{\max}(S_0^{-1}(BA_0^{-1}B^T)).
\end{align*}
Then the following holds for the \bpm{} and \bpp{} preconditioned matrices:
\begin{enumerate}[(i)]
    \item If $\cH_{1,1}$ is positive definite, then $\sigma(\cS_{1,1})\subseteq [-a_-,-b_-]\cup[c_-,d_-]$ with 
    \begin{align*}
        a_- &= \frac{\lambda^S_{\max}}{2} + \sqrt{\frac{(\lambda^S_{\max})^2}{4}  + \sigma_{\max} },
        & 
        b_- &= \lambda^S_{\min},
        & 
        c_- &=  \lambda^A_{\min},
        & 
        d_- &= 1.
    \end{align*}

    \item We have $\sigma(\cS_{-1,1})\subseteq [-a_+,-b_+]\cup[c_+,d_+]$ with 
    \begin{align*}
        a_+ &= \max\{1,\lambda^C_{\max}\}, &
        b_+ &= \frac{\lambda^S_{\min} }{(1 + \lambda^A_{\max})(1 +  \sigma_{\max})}, 
        \\
        c_+ &=  \lambda^A_{\min},&
        d_+ &= \frac{\lambda^T_{\max}}{2}+ \sqrt{\frac{(\lambda^T_{\max})^2}{4}+\sigma_{\max}} .
    \end{align*}
\end{enumerate}
\end{theorem}
\begin{proof}
We start with the preconditioned matrix $\cS_{\alpha,1}$ in \eqref{eqn:PinvS} for general $\alpha$.
Since $A_0 $ and $S_0$ are positive definite, this matrix  
is similar to 
    \begin{align*}
        \widetilde{\cS}_{\alpha,1} 
        &= 
        \begin{bmatrix}
            A_0^{-\frac12}AA_0^{-\frac12} &
            A_0^{-\frac12}B^TS_0^{-\frac12} \\ 
            S_0^{-\frac12}BA_0^{-\frac12}(I_n-\alpha A_0^{-\frac12}AA_0^{-\frac12}) & 
            -S_0^{-\frac12}C S_0^{-\frac12} -\alpha S_0^{-\frac12}BA_0^{-1}B^T S_0^{-\frac12}
        \end{bmatrix}
        \\
        &=
        \begin{bmatrix}
            \widetilde{A} & \widetilde{B}^T \\ 
            \widetilde{B}(I_n-\alpha \widetilde{A}) & 
            -(\widetilde{C} +\alpha \widetilde{B}\widetilde{B}^T)
        \end{bmatrix}.
    \end{align*}
    Thus, if $[x^T\;y^T]^T$ is an eigenvector of $\widetilde{\cS}_{\alpha,1}$ we get the eigenvalue equations 
    \begin{subequations}
    \begin{align}
        \lambda x &= \widetilde{A}x + \widetilde{B}^Ty,
        \label{eq:eig_BPv1}
        \\
        \lambda y &= \widetilde{B}(I_n-\alpha \widetilde{A})x - (\widetilde{C}+\alpha \widetilde{B}\widetilde{B}^T)y. \label{eq:eig_BPv2}
    \end{align}
    \end{subequations}
    Substituting  \eqref{eq:eig_BPv1} into \eqref{eq:eig_BPv2} for $\widetilde{B}^Ty$, we obtain 
    \begin{equation}\label{eq:eig_BPv3}
        \lambda y 
        = \widetilde{B}(I_n-\alpha \widetilde{A})x - \widetilde{C}y - \alpha \widetilde{B}(\lambda I -\widetilde{A})x
        = (1-\alpha\lambda)\widetilde{B}x -\widetilde{C}y.
    \end{equation}

Our first goal is to show that the positive eigenvalues of the \bpm{} and \bpp{} preconditioned matrices lie in the intervals $[c_{-},d_{-}]$ and $[c_+,d_+]$ stated in (i) and (ii), respectively. Thus, suppose that $\lambda>0$, then $(\lambda I_m +\widetilde{C})\succ0$, and we see from \eqref{eq:eig_BPv3} that 
    \begin{equation}\label{eq:eig_BPv_y}
    y = (1-\alpha\lambda)(\widetilde{C}+\lambda I_m)^{-1}\widetilde{B}x.
    \end{equation}
    In particular,  
    $x\neq 0$, because otherwise $x=0$ and $y=0$, which is a contradiction. 
    Substituting \eqref{eq:eig_BPv_y} into \eqref{eq:eig_BPv1}, yields 
    \begin{equation*}
        \lambda x = \wt{A}x + (1-\alpha\lambda)\wt{B}^T(\widetilde{C}+\lambda I_m)^{-1}\widetilde{B}x.
    \end{equation*}
 From this we obtain 
\begin{equation}\label{eq:eig_BPv_pos}
        \lambda = \frac{x^T(\wt{A}+ (1-\alpha\lambda)\wt{B}^T(\widetilde{C}+\lambda I_m)^{-1}\widetilde{B})x}{\norm{x}^2}.
    \end{equation}

For the \bpm{} preconditioned matrix, i.e., $\alpha=1$,     
we see from \eqref{eq:eig_BPv_pos} that 
    \begin{equation*}
        \lambda = \frac{x^T(\wt{A}+\wt{B}^T(\widetilde{C}+\lambda I_m)^{-1}\widetilde{B})x}{x^T(I_n+\wt{B}^T(\widetilde{C}+\lambda I_m)^{-1}\widetilde{B})x}.
    \end{equation*}
    Since $\cH_{1,1}$ is positive definite, we have $I_n-\widetilde{A}\succ0$, and hence
    \begin{equation*}
        \lambda \le \frac{x^T(I_n+\wt{B}^T(\widetilde{C}+\lambda I_m)^{-1}\widetilde{B})x}{x^T(I_n+\wt{B}^T(\widetilde{C}+\lambda I_m)^{-1}\widetilde{B})x} \le 1=d_{-}.
    \end{equation*}
    On the other hand, since $\lambda_{\min}(\wt{A})<1$ we get that
    \begin{align*}
        \lambda 
        &= \frac{x^T(\wt{A} + \wt{B}^T(\widetilde{C}+ \lambda I_m)^{-1}\widetilde{B})x}{x^T(I_n + \wt{B}^T(\widetilde{C}+ \lambda I_m)^{-1}\widetilde{B})x}
        \\
        &
        \ge 
        \lambda_{\min}(\wt{A}) \;
        \frac{x^T( I_n + \wt{B}^T(\widetilde{C}+ \lambda I_m)^{-1}\widetilde{B})x}{x^T(I_n + \wt{B}^T(\widetilde{C}+ \lambda I_m)^{-1}\widetilde{B})x}
        \\
        &=
        \lambda_{\min}(\wt{A})=\lambda^A_{\min}=c_{-}.
    \end{align*}
    
Similarly, for the \bpp{} preconditioned matrix, i.e., $\alpha=-1$, Equation 
\eqref{eq:eig_BPv_pos} yields 
    \begin{equation*}
        \lambda = \frac{x^T(\wt{A}+ (1+\lambda)\wt{B}^T(\widetilde{C}+\lambda I_m)^{-1}\widetilde{B})x}{\norm{x}^2}.
    \end{equation*}
    Since $\lambda>0$, it immediately follows    
    that $\lambda \ge \lambda_{\min}(\wt{A})=\lambda^A_{\min}=c_{+}$.
    To obtain an upper bound on the positive eigenvalues, we estimate $(\wt{C}+\lambda I)^{-1}\le \frac{1}{\lambda} I_m $   to obtain 
    \begin{equation*}
        \lambda\le \frac{x^T\wt{A}x}{\norm{x}^2} + \frac{\lambda+1}{\lambda} \frac{x^T\wt{B}^T\wt{B}x}{\norm{x}^2}.
    \end{equation*}
    Therefore, 
    \begin{equation*}
        \lambda^2 - \lambda\; \lambda_{\max}(\wt{A}+\wt{B}^T\wt{B}) - \lambda_{\max}(\wt{B}^T\wt{B})\le 0,
    \end{equation*}
    which yields 
    \begin{equation*}
        \lambda \le \frac{\lambda_{\max}(\wt{A}+\wt{B}^T\wt{B})}{2}+ \sqrt{\frac{\lambda_{\max}(\wt{A}+\wt{B}^T\wt{B})^2}{4}+\lambda_{\max}(\wt{B}^T\wt{B})}=d_+.
    \end{equation*}

Next, we show that the negative eigenvalues of the \bpm{} and \bpp{} preconditioned matrices lie in the intervals $[-a_{-},-b_{-}]$ and $[-a_+,-b_+]$ stated in (i) and (ii), respectively.
    To that end, assume $\lambda<0$. In this case $\wt{A}-\lambda I_n\succ0$ and thus using \eqref{eq:eig_BPv1} we have 
    \begin{equation}\label{eq:eig_BPv4}
        x = - (\wt{A}-\lambda I_n)^{-1}\wt{B}^Ty.
    \end{equation}
    We must have $y\neq 0$, because otherwise we would have $x=0$, which leads to a contradiction. Substituting \eqref{eq:eig_BPv4} into \eqref{eq:eig_BPv3} and multiplying on the left by $y^T$ leads to
\begin{equation}\label{eq:eig_BBv_neg}
        -\lambda = \frac{y^T(\wt{C}+(1-\alpha\lambda)\wt{B}(\wt{A}-\lambda I_n)^{-1}\wt{B}^T)y}{\norm{y}^2}.
    \end{equation}
    
    For $\alpha=1$, we can bound \eqref{eq:eig_BBv_neg} from below by using  $I_n\succ\wt{A}$, so that 
    \[
    -\lambda \ge  \frac{y^T(\wt{C}+\wt{B}\wt{B}^T)y}{\norm{y}^2} \ge \lambda_{\min}(\wt{C}+\wt{B}\wt{B}^T)
    =\lambda^S_{\min}=b_{-}.
    \]
    To obtain an upper bound, we use the estimate 
    $(\wt{A}-\lambda I_n)^{-1}\le -\frac{1}{\lambda} I_n $, and obtain 
    \[
    -\lambda \le \frac{y^T\wt{C}y}{\norm{y}^2} + \frac{1+(-\lambda)}{(-\lambda)}\frac{y^T\wt{B}\wt{B}^Ty}{\norm{y}^2}.
    \]
    Similarly as before, we see that
    \[
    -\lambda \le \frac{\lambda_{\max}(\wt{C}+\wt{B}\wt{B}^T)}{2} + \sqrt{\frac{\lambda_{\max}(\wt{C}+\wt{B}\wt{B}^T)^2}{4}  + \lambda_{\max}(\wt{B}\wt{B}^T)}
    =a_{-}.
    \]

    For $\alpha=-1$, we use \eqref{eq:eig_BBv_neg} together with $\wt{C}\le \lambda_{\max}(\wt{C}) I_m$ to obtain 
    \[
    -\lambda = \frac{y^T(\wt{C}+\wt{B}(\wt{A}-\lambda I_n)^{-1}\wt{B}^T)y}{y^T(I_m+\wt{B}(\wt{A}-\lambda I_n)^{-1}\wt{B}^T)y}
    \le \max\{1,\lambda_{\max}(\wt{C})\}
    =a_{+}.
    \]
    Obtaining a lower bound on $-\lambda$ is  trickier in this case, because $\wt{C}$ is only 
    semi-definite. Note that by \eqref{eq:eig_BBv_neg}, we have 
    \begin{equation*}
        y^T\wt{C}y = y^T(-\lambda I_m - (1+\lambda )\wt{B}(\wt{A}-\lambda I_n)^{-1}\wt{B}^T) y.
    \end{equation*}
    Thus,  
    \begin{align*}
        -\lambda - & \frac{y^T(\wt{C}+\wt{B}(\wt{A}+ I_n)^{-1}\wt{B}^T)y}{y^T(I_m+\wt{B}(\wt{A}+ I_n)^{-1}\wt{B}^T)y}
        \\
        & = 
        -\lambda - \frac{y^T(-\lambda I_m - (1+\lambda )\wt{B}(\wt{A}-\lambda I_n)^{-1}\wt{B}^T +\wt{B}(\wt{A}+ I_n)^{-1}\wt{B}^T)y}{y^T(I_m+\wt{B}(\wt{A}+ I_n)^{-1}\wt{B}^T)y}
        \\
        & = \frac{ (-\lambda-1 )y^T\wt{B}((\wt{A}+ I_n)^{-1} -(\wt{A}- \lambda I_n)^{-1})\wt{B}^Ty}
        {y^T(I_m+\wt{B}(\wt{A}+ I_n)^{-1}\wt{B}^T)y}
        \ge0,
    \end{align*}
    since $-\lambda \ge 1$ if and only if $(\wt{A}+ I_n)^{-1} \ge (\wt{A}-\lambda I_n)^{-1}$.
    Using  $0\prec \wt{A}$, we finally obtain 
    \begin{equation*}
        -\lambda \ge 
        \frac{y^T(\wt{C}+\wt{B}(\wt{A}+ I_n)^{-1}\wt{B}^T)y}{y^T(I_m+\wt{B}(\wt{A}+ I_n)^{-1}\wt{B}^T)y}
        \ge
        \frac{\lambda_{\min}(\wt{C}+\wt{B}\wt{B}^T)}{(1+\lambda_{\max}(\wt{A})) (1+\lambda_{\max}(\wt{B}\wt{B}^T))}=b_+, 
    \end{equation*}
    which completes the proof.
\end{proof}

Similarly to \cref{thm:eig_BP}, the bounds in \cref{thm:eig_bounds} depend on the extreme eigenvalues of $A_0^{-1}A$ and $S_0^{-1}(C+BA_0^{-1}B^T)$. Moreover, the bounds get tighter when these matrices are close to the identity matrix.

\cref{thm:eig_bounds} shows that the positive eigenvalues of the \bpm{} preconditioned matrix $\cS_{1,1}$ can be bounded only in terms of the extreme eigenvalues of $A_0^{-1}A$. Hence, when $A_0$ is a good approximation to $A$, then the positive eigenvalues of $\cS_{1,1}$ will be tightly clustered, approaching $1$ from below. We already made a similar observation in the simple one-parameter problem above.

\begin{example}\label[example]{ex:special-case}
    We apply the bounds from \cref{thm:eig_BP,thm:eig_bounds} to the special case $C=0$, $A_0 = \gamma A\succ0$ and $S_0=BA^{-1}B^T\succ0$. Doing so, we obtain 
    \begin{equation*}
    \begin{aligned}
    {\rm \bp{}}:&& \sigma(\cS_{1,-1}) &\subseteq 
        \left[ \frac{\gamma^2}{2(1+\gamma)}, \frac{2}{\gamma} \right]
        \quad \text{for } 0<\gamma <1,
        \\
        {\rm BP}^{+}:&& \sigma(\cS_{-1,1}) &\subseteq 
        \left[ -1, -\frac{1}{2+\gamma+\frac{1}{\gamma}}  \right]
        \cup
        \left[ \frac{1}{\gamma}, \frac{1}{\gamma}(1 + \sqrt{1+\gamma}) \right]
        \quad \text{for } \gamma >0,\\
        {\rm BP}^{-}:&&\sigma(\cS_{1,1}) &\subseteq 
        \left[ -\frac{1}{2\gamma}(1+ \sqrt{1 + 4\gamma}), -\frac{1}{\gamma} \right]
        \cup
        \left[ \frac{1}{\gamma}, 1 \right]
        \quad \text{for } \gamma >1.
    \end{aligned}
    \end{equation*}
    Interestingly, in this case the bounds from \cref{thm:eig_bounds} on the positive eigenvalues of $\cS_{-1,1}$ are tight; see \eqref{eq:eigvals}. 
    In Figure~\ref{fig:bounds-special-case} we plot the eigenvalues and the corresponding bounds from \cref{thm:eig_BP,thm:eig_bounds} for the special case considered here and different ranges of $\gamma>0$. The figure shows that the bounds describe the  behavior  of the eigenvalues rather faithfully, in particular the ones in \cref{thm:eig_bounds} shown in Figures~\ref{fig:bounds-special-case}~(b) and (c).
\begin{figure}
    \centering
\hfill
\begin{subfigure}{.28\linewidth}
    \hspace{-3em}
    \includegraphics{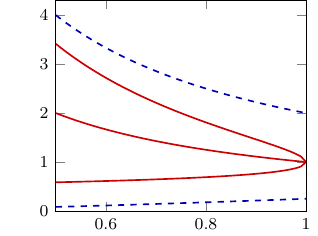}    
    \caption{BP}
\end{subfigure}
\hfill
\begin{subfigure}{.28\linewidth}
    \hspace{-3em}
    \includegraphics{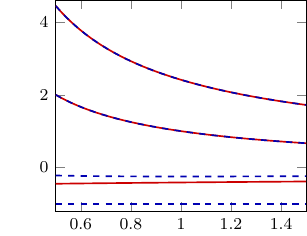}    
    \caption{\bpp}
\end{subfigure}
\hfill
\begin{subfigure}{.28\linewidth}
    \hspace{-3em}
    \includegraphics{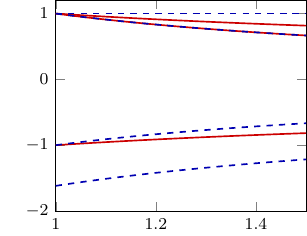}    
    \caption{\bpm}
\end{subfigure}
\hfill

    \caption{Eigenvalues (solid red) and bounds (dashed blue) from \cref{thm:eig_BP,thm:eig_bounds} for  
    \cref{ex:special-case}   
    with different values of $\gamma$.} 
\label{fig:bounds-special-case}  
\end{figure}
\end{example}

\section{Examples} 
\label{sec:examples}

In this section, we study the effectiveness of the three different BP-type preconditioning methods, also compared to common block-diagonal preconditioning. All numerical results were computed using \textsc{Matlab} version R2026a on an Apple M5 with 16 GB of unified memory.

For our experiments, we have used our own \textsc{Matlab} implementations of the $H$-PCG and $H$-PMINRES methods provided in Appendix \ref{sec:methods}. All experiments used $x_0=0$ as the initial vector in iterative methods.
In our implementation of $H$-PMINRES the $H$-inner products are computed without explicitly forming $A_0$ and $S_0$, as explained in \cite[Remark~7]{BraPas1988} and \cite[Section~3.2]{Sto2009}. Although this is not necessary in \cref{ex:exact_stokes,ex:gauss-seidel}, where $A_0$ is explicitly available,  in \cref{ex:multigrid}  we  have the action of $A_0$ only implicitly available through $A_0^{-1}$.
Recall that the $H$-PCG method  
requires that the preconditioned matrix be (positive)
definite in the $H$-inner product, which is the case for the original BP method. The $H$-PMINRES method, on the other hand, is applicable when the preconditioned matrix is only symmetric (possibly indefinite) in the $H$-inner product, which holds for all three BP-type methods. Each of the three techniques is based on a different inner product, and hence the corresponding iterative methods are optimal in different norms. 
In addition, we compare the three methods to common block-diagonal preconditioning, i.e., $\cP=\diag(A_0,S_0)$, 
where we apply the preconditioned MINRES method (also through our own implementation). 
Note that since $\cP=\cP^T\succ0$, the matrix $\cP^{-1}\cS$ is $\cP$-selfadjoint, 
and hence preconditioned MINRES is optimal in the $\cP$-norm (i.e., the preconditioned residual is minimized in this norm).
In order to provide a thorough comparison of all methods, we present the computed relative Euclidean norms ($\|b-\cS x_k\|_2/\|b\|_2$) as well as the  computed \emph{preconditioned} relative $H$-norms 
($\|\cP^{-1}b-\cP^{-1}\cS x_k\|_H/\|\cP^{-1}b\|_H$).

The saddle point systems we consider are generated by finite element discretizations of the Stokes model problem 
\begin{equation*}
\begin{aligned}
    - \nabla^2 u + \nabla p &= f, & \text{in } \Omega,
    \\
    \nabla \cdot u &= 0 , & \text{in } \Omega,
\end{aligned}
\end{equation*}
with different boundary conditions as described in \cite[Section 3]{ElmSilWat2014}. The discretizations of all the problems considered were done using IFISS \cite{ifiss} and resulted in saddle point systems
of the form 
\begin{equation}\label{eq:stokes-sys}
        \begin{bmatrix} A & B^T \\ B & -C \end{bmatrix}
        \begin{bmatrix}u \\ p \end{bmatrix}
        =\begin{bmatrix}f \\ g \end{bmatrix}.
\end{equation}
Throughout all experiments, we choose the block $S_0$ as the so-called pressure mass matrix, also generated by IFISS, which in this application is diagonal~\cite[p.~198]{ElmSilWat2014} and known to be spectrally equivalent to the Schur complement $S$; see, e.g.,~\cite[Section~4.2.1]{ElmSilWat2014}.

\begin{example}\label[example]{ex:exact_stokes}
    We consider the steady horizontal flow in a channel  \cite[Example~3.1.1]{ElmSilWat2014} and use IFISS default parameters\footnote{
    Channel domain: $1$, grid parameter: $4$, $Q_1$-$P_0$-elements, stabilization parameter: $1/4$, uniform streamlines. 
    }
    for the discretization. The resulting saddle point matrix in 
    \eqref{eq:stokes-sys} is of size $834 \times 834$ with $n = 578 $ and $m = 256$. In this experiment we set $A_0 = \gamma A$, where $\gamma$ is chosen such that the respective inner product matrices are positive definite for each of the three BP-type methods.   

    Figures~\ref{fig:exact_stokes_BP} (a) to (c) 
    and Figures~\ref{fig:exact_stokes_BP_Hnorm} (a) to (c)
    show the computed relative Euclidean norms ($\|b-\cS x_k\|_2/\|b\|_2$) as well as the computed preconditioned relative $H$-norms ($\|\cP^{-1}b-\cP^{-1}\cS x_k\|_H/\|\cP^{-1}b\|_H$), respectively, for all BP-type methods and different parameters $\gamma$. 
    Each figure contains a comparison with block-diagonal preconditioning (and the MINRES method).
    
    Our first observation is that for each method the convergence curves for the Euclidean norms and the preconditioned $H$-norms look very similar. One noticeable but overall small difference is that curves for the $H$-norms are monotonically decreasing (unlike for the Euclidean norms), since the $H$-norms are minimized in each step by 
    $H$-PMINRES.
    Moreover, we observe that the convergence behavior for the \bp{} and \bpm{} preconditioners depends strongly on the choice of $\gamma$, i.e., how well $A_0$ approximates $A$, while the performance for the \bpp{} preconditioner is much less sensitive to this choice. In fact, for all choices of $\gamma$ the performance for the \bpp{} preconditioner is similar to the one for block-diagonal preconditioning. This behavior is consistent with our observations in \cref{sec:special-case} about the structural similarity of the spectra of the corresponding preconditioned matrices $\cS_{-1,1}$ and $\cS_{0,1}$.
    
    In Figure~\ref{fig:exact_stokes_BP} (d) and Figure~\ref{fig:exact_stokes_BP_Hnorm} (d), the different preconditioning strategies are compared against each other for respective choices of $\gamma$ that lead to well defined methods and $A_0\approx A$, namely $\gamma=0.9$ for \bp{}, $\gamma = 1$ for \bpp{}, and $\gamma = 1.1$ for \bpm{}. 
    For the \bp{} preconditioner, we show the convergence curves of $H$-PCG and $H$-PMINRES (which are almost indistinguishable).
    The figures indicate that the positive definite system obtained by the BP preconditioner seems to be advantageous to the indefinite system obtained by the \bpm{} preconditioner, but that both outperform the \bpp{} preconditioner (and block-diagonal preconditioning).

    In Figure~\ref{fig:exact_stokes_eig_BP} we show the numerically computed spectra (obtained using \textsc{Matlab}'s \texttt{eig} function) of the preconditioned saddle point matrices for block-diagonal preconditioning and three BP-type methods with the parameter choices used in Figure~\ref{fig:exact_stokes_BP} (d) and Figure~\ref{fig:exact_stokes_BP_Hnorm} (d). The dashed lines in Figures~\ref{fig:exact_stokes_eig_BP} (a) to (c) indicate the eigenvalue bounds for the BP-type methods from \cref{thm:eig_BP} and \cref{thm:eig_bounds}.
    The spectra for the block-diagonal preconditioner and for the \bpp{} preconditioner are of a similar structure (cf.~our remarks above), which may also explain the rather similar convergence behavior. On the other hand, as discussed after \cref{thm:eig_bounds}, the positive eigenvalues of the \bpm{} preconditioner are tightly clustered below $1$.     
\begin{figure}[]
\centering

\hfill
\begin{subfigure}[t]{0.49\textwidth}
\centering  
\hspace{-1.5ex}
\includegraphics{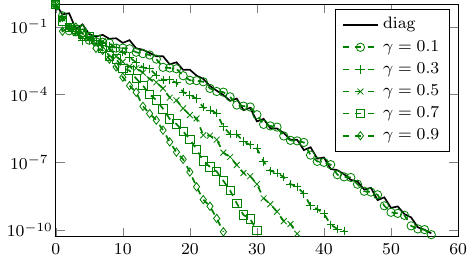}
    \caption{BP} 
    \label{sfig:exact_stokes_BP}
\end{subfigure}
\hfill
\begin{subfigure}[t]{0.49\textwidth}
\centering
\hspace{-1.5ex}
\includegraphics{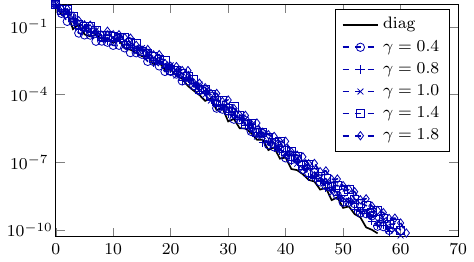}
    \caption{\bpp} 
    \label{sfig:exact_stokes_BPp}
\end{subfigure}
\hfill

\hfill
\begin{subfigure}[t]{0.49\textwidth}
         \centering
         \hspace{-1.5ex}
\includegraphics{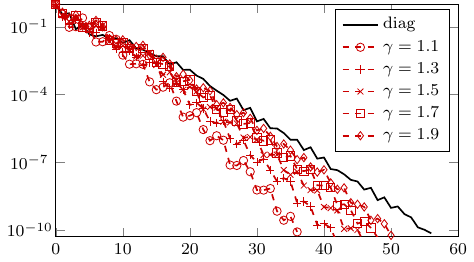}
    \caption{\bpm}
    \label{sfig:exact_stokes_BPm}
\end{subfigure}
\hfill
\begin{subfigure}[t]{0.49\textwidth}
         \centering
         \hspace{-1.5ex}
\includegraphics{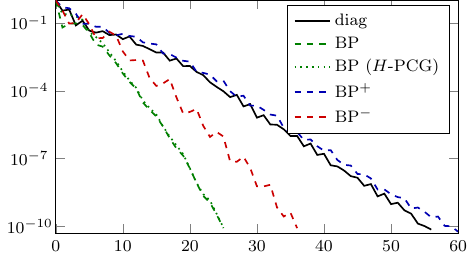}
    \caption{Comparison of the different BP methods with $\gamma=0.9$ for BP, $\gamma=1$ for \bpp{} and $\gamma=1.1$ for \bpm. }
    \label{sfig:exact_stokes_comp}
\end{subfigure}
\hfill
    
\caption{Comparison of the (unpreconditioned) relative Euclidean residual norms $\|b-\cS x_k\|_2/\|b\|_2$ in \cref{ex:exact_stokes} for different 
BP-type preconditioners and a block-diagonal preconditioner in $H$-PMINRES (and $H$-PCG for \bp{} in (d)).}
\label{fig:exact_stokes_BP}
\end{figure}

\begin{figure}[]
\centering

\hfill
\begin{subfigure}[t]{0.49\textwidth}
         \centering
         \hspace{-1.5ex}
\includegraphics{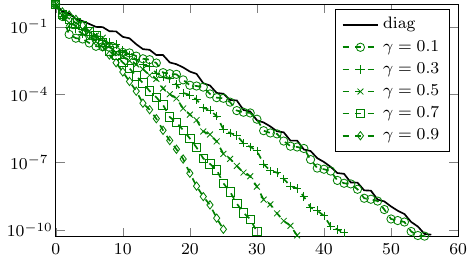}
    \caption{BP} 
    \label{sfig:exact_stokes_BPH}
\end{subfigure}
\hfill
\begin{subfigure}[t]{0.49\textwidth}
         \centering
         \hspace{-1.5ex}
\includegraphics{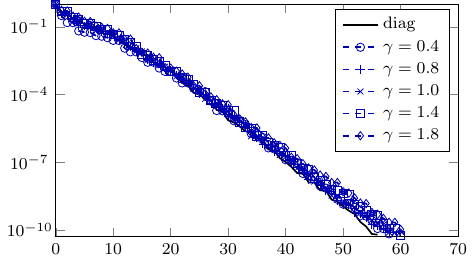}
    \caption{\bpp} 
    \label{sfig:exact_stokes_BPpH}
\end{subfigure}
\hfill

\hfill
\begin{subfigure}[t]{0.49\textwidth}
         \centering
         \hspace{-1.5ex}
    \includegraphics{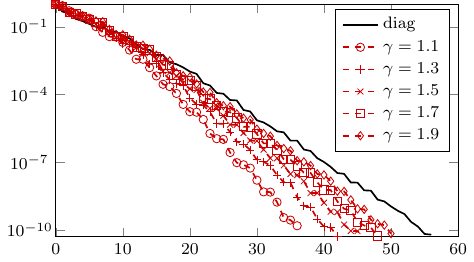}
    \caption{\bpm}
    \label{sfig:exact_stokes_BPmH}
\end{subfigure}
\hfill
\begin{subfigure}[t]{0.49\textwidth}
         \centering  
         \hspace{-1.5ex}
         \includegraphics{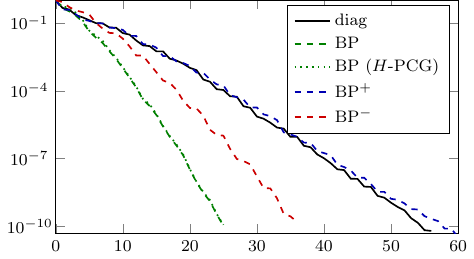}
    \caption{Comparison of the different BP-technique with $\gamma=0.9$ for BP, $\gamma=1$ for \bpp{} and $\gamma=1.1$ for \bpm. }
    \label{sfig:exact_stokes_compH}
\end{subfigure}
\hfill
    
\caption{Comparison of the  relative residual $H$-norms $\|\cP^{-1} b- \cP^{-1}\cS x_k\|_H/\|\cP^{-1}b\|_H$ in \cref{ex:exact_stokes} for different 
BP-type preconditioners and a block-diagonal preconditioner in $H$-PMINRES (and $H$-PCG for \bp{} in (d)).}
\label{fig:exact_stokes_BP_Hnorm}
\end{figure}

\begin{figure}[h]
\centering
\hfill
\begin{subfigure}[t]{0.49\textwidth}
         \centering
    \hspace{-2em}
    \includegraphics{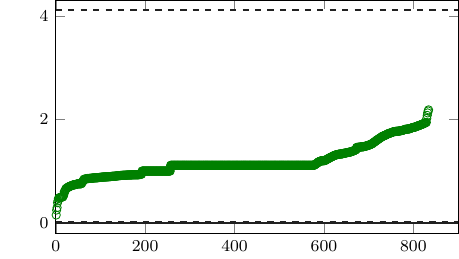}
    \caption{BP} 
    \label{sfig:exact_stokes_eig_BP}
\end{subfigure}
\hfill
\begin{subfigure}[t]{0.49\textwidth}
         \centering
    \hspace{-2em}
    \includegraphics{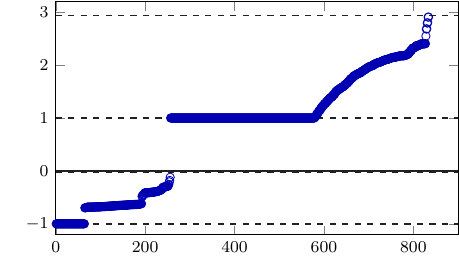}
    \caption{\bpp}
    \label{sfig:exact_stokes_eig_BPp}
\end{subfigure}
\hfill

\hfill
\begin{subfigure}[t]{0.49\textwidth}
         \centering
    \hspace{-2em}
    \includegraphics{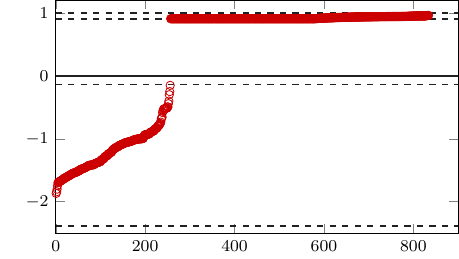}
    \caption{\bpm }
    \label{sfig:exact_stokes_eig_BPm}
\end{subfigure}
\hfill
\begin{subfigure}[t]{0.49\textwidth}
    \centering
    \hspace{-2em}
    \includegraphics{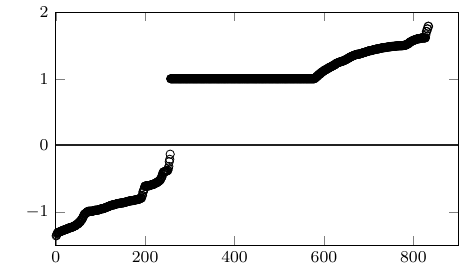}
    \caption{block-diagonal} 
    \label{sfig:exact_stokes_eig_diag}
\end{subfigure}
\hfill
    
\caption{Plots (a)-(c) show the eigenvalues of the preconditioned saddle point matrices for the three different BP-type preconditioners  in \cref{ex:exact_stokes} with the parameters $\gamma=0.9$ (\bp{}), $\gamma=1$ (\bpp{}), and $\gamma=1.1$ (\bpm{}) as used in Figure~\ref{fig:exact_stokes_BP}~(d) and Figure~\ref{fig:exact_stokes_BP_Hnorm}~(d). The dashed lines in plots (a)-(c) show the eigenvalue bounds from \cref{thm:eig_BP,thm:eig_bounds}. Plot (d) shows the eigenvalues of the block-diagonally preconditioned matrix (without bounds). 
}
\label{fig:exact_stokes_eig_BP}
\end{figure}

\end{example}

In practical computations using a scalar multiple $A_0=\gamma A$  
as a preconditioner for $A$ 
is too expensive. Thus, we now consider two examples where other choices  $A_0\approx A$ are used. 
Recall that in order to apply the $H$-PMINRES method to preconditioned systems, we need to ensure that the preconditioner $A_0$ satisfies $A_0\prec A$ for BP or $A_0\succ A$ for \bpm. 
Thus, in both cases the preconditioner should be chosen carefully. 
Satisfying $A_0\prec A$ is, in general, more difficult and usually requires rather precise knowledge of the spectra of $A$ and $A_0$;
see, e.g., \cite{MueUllLan2019}. The following two examples illustrate that there exist 
cases where $A_0\succ A$ is guaranteed  
without the need for additional spectral information. 

\begin{figure}[]
\centering
\hfill
\begin{subfigure}[t]{0.49\textwidth}
\centering
\hspace{-1.5ex}
\includegraphics{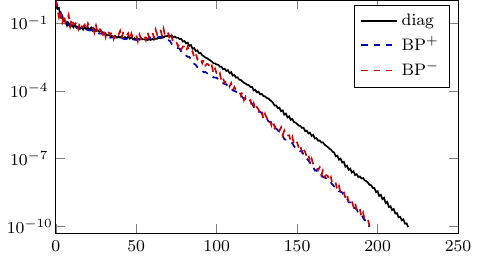}
\caption{$\|b-\cS x_k\|_2/\|b\|_2$}   
\label{sfig:GS_stokes_2norm}
\end{subfigure}
\hfill
\begin{subfigure}[t]{0.49\textwidth}
\centering
\hspace{-1.5ex}
\includegraphics{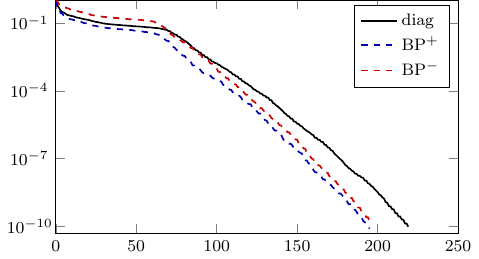}
\caption{$\|\cP^{-1}b-\cP^{-1}\cS x_k\|_H/\|\cP^{-1}b\|_H$} 
\label{sfig:GS_stokes_Hnorm}
\end{subfigure}
\hfill

\caption{Comparison of \bpp{}, \bpm{} and block-diagonal preconditioning   described in \cref{ex:gauss-seidel}.} 
\label{fig:GS_stokes}
\end{figure}

\begin{example}\label[example]{ex:gauss-seidel}
    We consider the flow over a backward facing step \cite[Example~3.1.2]{ElmSilWat2014} where we use the default parameters\footnote{
    horizontal dimensions: $[-1,5]$, grid parameter: $4$, grid stretch factor: 1, $Q_1$-$P_0$-elements, stabilization parameter: $1/4$.
    }
    in IFISS. The resulting saddle point matrix in 
    \eqref{eq:stokes-sys} is  of size $2,242 \times 2,242$ with $n=1,538$ and $m=704$. 
    In this experiment 
    we consider $A_0^{-1}$ to be one iteration of a scaled symmetric Gauss-Seidel method~\cite[Section 4.1]{Saa2003}.   
    For $A = L + D + L^T$, where $L$ is the strictly lower triangular part of $A$ and $D= \diag(A)$ is the diagonal part of $A$, we set 
    \begin{equation*}
        A_0 = \varphi  (D+L)D^{-1}(D+L^T) =  \varphi A + \varphi LD^{-1}L^T.
    \end{equation*}
     Now $A\succ0$ implies   
     $D \succ0$ and thus $A_0\succ A$ is automatically satisfied for any $\varphi>1$, so that the \bpm{} method is applicable. For our computations we choose $\varphi=1.1$, and take the same $A_0$ also for the block-diagonal and the \bpp{} preconditioner. In order to use the \bp{} method we have to guarantee $A\succ A_0$, which would require an appropriate scaling based on additional spectral information. We therefore do not consider the \bp{} method in this example.
    Figure~\ref{fig:GS_stokes} compares the relative residual norms obtained from $H$-PMINRES for the \bpp, \bpm{} and block-diagonal preconditioners. 
    While all methods require a substantial number of iterations, the \bpp{} and \bpm{} methods exhibit similar convergence behavior, while both lead to  fewer iterations than with block-diagonal preconditioning. Computing the extreme eigenvalues of $A_0^{-1}A$ using \textsc{Matlab}'s \texttt{eig} function reveals that 
    \begin{equation*}
    \lambda_{\min}(A_0^{-1}A) \approx 0.0519
    \quad\text{and}\quad 
    \lambda_{\max}(A_0^{-1}A) \approx 0.9091,
    \end{equation*}
    indicating that $A_0$ is not a close 
    approximation to $A$.  
    We have already seen in \cref{ex:exact_stokes} that for a poor choice of $A_0$ the \bpp{} and \bpm{} methods result in a similarly slow convergence behavior.
\end{example}

\begin{example}\label[example]{ex:multigrid} 
    We again consider the flow over a backward facing step like in \cref{ex:gauss-seidel}, now with the IFISS grid parameter $6$. The resulting  system is of size $34,306\times 34,306$ with $n=23,042$ and $m=11,264$. 

    The matrix $A$ in \eqref{eq:stokes-sys} is of the form $A=\diag(\wt{A},\wt{A})$, where $\wt{A}$ is a scalar Laplace matrix; see \cite[p.~130]{ElmSilWat2014}. For the (1,1)-block of the preconditioner we use the matrix 
    \begin{equation*}
        A_0^{-1} = \frac{1}{\varphi}\diag(M^{-1},M^{-1}),
    \end{equation*}
    where $M^{-1}$ is the action of a single multigrid V-cycle with two sweeps of symmetric Gauss-Seidel smoothing, and an exact solve on the coarsest level. 
    In our implementation, the multigrid V-cycle is realized by using the function {\rm \texttt{amg\_v\_cycle}} provided by IFISS. 
    The required grids are generated via {\rm \texttt{amg\_grids\_setup2D($\wt{A}$)}}, and the parameters by using the function {\rm \texttt{amg\_smoother\_params(grid,'PGS',2)}}.
    
    As explained in \cite[Section~5.1]{Vas2008}, the operator $M$ is symmetric and satisfies $\wt{A}\preceq M$, so that the condition $A_0-A\succ0$ required for \bpm{} is satisfied for any $\varphi>1$. As in \cref{ex:gauss-seidel} we use $\varphi = 1.1$. (Note that $M-\wt{A}$ is in general only positive semi-definite and hence the additional scaling is necessary.)

    Figure~\ref{fig:amg_stokes} 
    again compares the relative residual norms obtained from $H$-PMINRES for the \bpp, \bpm{} and block-diagonal preconditioners. We observe that the $A_0$ in this example gives
    a significant improvement over the Gauss-Seidel preconditioner from \cref{ex:gauss-seidel}. Note in particular that here we consider a significantly larger linear algebraic system.  
    Moreover, the \bpm{} preconditioner outperforms the two other techniques, which may be attributed to the fact that $A_0$ is a very good preconditioner for $A$.

\begin{figure}[]
\centering
\hfill
\begin{subfigure}[t]{0.49\textwidth}
\centering
\hspace{-1.5ex}
\includegraphics{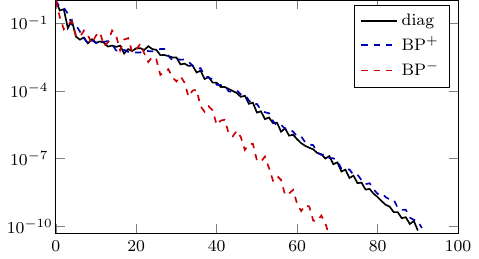}
\caption{$\|b-\cS x_k\|_2/\|b\|_2$}     
\label{sfig:amg_stokes_2norm}
\end{subfigure}
\hfill
\begin{subfigure}[t]{0.49\textwidth}
\centering
\hspace{-1.5ex}
\includegraphics{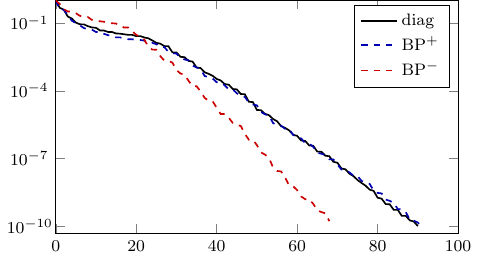}
\caption{$\|\cP^{-1}b-\cP^{-1}\cS x_k\|_H/\|\cP^{-1}b\|_H$} 
\label{sfig:amg_stokes_Hnorm}
\end{subfigure}
\hfill

\caption{Comparison of \bpp{}, \bpm{} and block-diagonal preconditioning  as described in \cref{ex:multigrid}.} 
\label{fig:amg_stokes}
\end{figure}

\end{example}

\section{Concluding remarks} 
\label{sec:conclusion}

We have analyzed and compared three different BP-type preconditioners, namely the original \bp{} method~\cite{BraPas1988}, \bpp{}~\cite{StoWat2008}, and \bpm{}. The three methods differ by sign choices in their inner product matrix and corresponding choices in the preconditioner itself. In each case, the nonsymmetric preconditioned matrix is selfadjoint with respect to the given non-standard inner product. Hence, this matrix is diagonalizable with real eigenvalues, which allows us to apply optimal Krylov subspace methods with efficient three-term recurrences to solve the preconditioned linear algebraic systems.

Although our analysis does not result in a strict recipe for the choice of a BP-type preconditioner, we have made some observations that may guide such a choice. The original \bp{} method is the most restrictive, but leads to a positive definite preconditioned matrix that allows the application of a generalized conjugate gradient method ($H$-PCG). Thus, if the original \bp{} method is applicable, it can be quite effective. In contrast, the preconditioned matrices in the \bpm{} and \bpp{} methods remain indefinite. The \bpp{} method has the fewest restrictions, but structurally is comparable to common block-diagonal preconditioning and did not outperform the latter in our numerical experiments. The \bpm{} method, which so far has not been studied in detail in the literature, has several interesting mathematical properties and in practical applications appears to be an attractive alternative to the original \bp{} method due to its less restrictive requirements. Note, in particular, that the necessary condition $A_0\succ A$ can be achieved even for singular or indefinite matrices $A$, which we have not studied in detail in this paper.

As shown in \cref{lem:basic}, several other variants of the three basic techniques can be constructed by considering general parameters $\alpha,\beta\in\R$ (and $\beta\neq 0$) instead of just $\alpha,\beta\in\{-1,1\}$. Yet, we believe that the main features of the class of BP-type preconditioners are represented by the three basic techniques, and that these features are largely revealed by our current analysis. Thus, unless there is a compelling practical reason for a specific choice of $\alpha$ and $\beta$ in a given application, it is likely sufficient to consider the three basic techniques.

\section*{Acknowledgments}
We thank Ludwig Rooch for his helpful comments regarding multigrid. ChatGPT was used for language polishing and help with coding the numerical examples and with drafting \cref{thm:eig_bounds}.
The authors assume responsibility for all content.
Part of this work was supported by the Simons Institute for the Theory of Computing, and conducted when the first and second authors were visiting the Institute in Fall 2025.



\pagebreak
\appendix

\crefalias{section}{appendix}

\section{Optimal three-term recurrence methods in the $H$-inner product}\label{sec:methods} 

\medskip
\Cref{thm:optimal} gives necessary and sufficient conditions for the existence of an optimal three-term recurrence method for a given matrix in a given inner product. For completeness we state two such methods, $H$-PCG from \cite{DolGouStoWat2010} (Algorithm \ref{alg:hpcg}) and $H$-PMINRES from \cite{PesWat2013} (Algorithm \ref{alg:hpminres}), which we have also used in the numerical experiments in this paper. We assume that the preconditioned matrix $P^{-1}M$ is $H$-selfadjoint, hence in particular $H$-normal(1). For $H$-PCG it is additionally assumed that $HP^{-1}M$ is positive definite.

\begin{algorithm}
\caption{$H$-PCG as in \cite[Algorithm 3.1]{DolGouStoWat2010}.}
\label{alg:hpcg}
\begin{algorithmic}
    \REQUIRE Nonsingular $M,P,H\in\C^{n\times n}$ with $H$ and
    $H P^{-1}M$ Hermitian positive definite, $b,x_0\in\C^{n}$.

    \STATE \textbf{Initialize:} Set $r_0 = P^{-1}(b-Mx_0)$
    and $p_0 = r_0$.

    \FOR{$k = 0, 1, \dots$}
        \STATE $\alpha_k =
        \langle r_k, p_k \rangle_H
        \,/\, \langle P^{-1} M p_k, p_k \rangle_H$
        \STATE $x_{k+1} = x_k + \alpha_k p_k$
        \STATE $r_{k+1} = r_k - \alpha_k P^{-1} M p_k$
        \STATE $\beta_{k+1} =
        \langle P^{-1} M r_{k+1}, p_k \rangle_H
        \,/\, \langle P^{-1} M p_k, p_k \rangle_H$
        \STATE $p_{k+1} = r_{k+1} - \beta_{k+1} p_k$
    \ENDFOR
\end{algorithmic}
\end{algorithm}

\begin{algorithm}
\caption{$H$-PMINRES as in \cite[Algorithm 2]{PesWat2013}}
\label{alg:hpminres}
\begin{algorithmic}
    \REQUIRE Nonsingular $M,P,H\in\C^{n\times n}$ with $H$
    Hermitian positive definite and $H P^{-1}M$ Hermitian,
    $b,x_0\in\C^{n}$.

    \STATE \textbf{Initialize:} Set
    $v_1 = b - Mx_0$,
    $z_1 = P^{-1}v_1$,
    $\gamma_1 = \norm{z_1}_{H}$,
    $v_0 = 0$, $w_0 = 0$, $w_1 = 0$,
    $\eta_1 = \gamma_1$,
    $s_0 = s_1 = 0$ and $c_0 = c_1 = 1$.

    \FOR{$k = 1, 2, \dots$}
        \STATE $v_k = v_k/\gamma_k$
        \STATE $z_k = z_k/\gamma_k$
        \STATE $\delta_k = \langle P^{-1}Mz_k, z_k \rangle_H$
        \STATE $v_{k+1} = Mz_k - \delta_k v_k - \gamma_k v_{k-1}$
        \STATE $z_{k+1} = P^{-1}v_{k+1}$
        \STATE $\gamma_{k+1} = \|z_{k+1}\|_H$
        \STATE $\alpha_0 = c_k \delta_k - c_{k-1} s_k \gamma_k$
        \STATE $\alpha_1 = \sqrt{\alpha_0^2 + \gamma_{k+1}^2}$
        \STATE $\alpha_2 = s_k \delta_k + c_{k-1} c_k \gamma_k$
        \STATE $\alpha_3 = s_{k-1} \gamma_k$
        \STATE $c_{k+1} = \alpha_0/\alpha_1,\quad
        s_{k+1} = \gamma_{k+1}/\alpha_1$
        \STATE $w_{k+1} =
        (z_k - \alpha_3 w_{k-1} - \alpha_2 w_k)/\alpha_1$
        \STATE $x_k = x_{k-1} + c_{k+1}\eta_k w_{k+1}$
        \STATE $\eta_{k+1} = -s_{k+1}\eta_k$
    \ENDFOR
\end{algorithmic}
\end{algorithm}

\end{document}